\documentclass[MSNbibl,number,citesort,seceqn,dvips]{arxbj}

\aid{0}
\volume{19}
\issue{5B}
\pubyear{2013}
\firstpage{2153}
\lastpage{2166}
\doi{10.3150/12-BEJ447} 

\makeatletter
\newcommand{\rrVert}{\Vert}
\newcommand{\llVert}{\Vert}
\newcommand{\eqref}[1]{(\ref{#1})}
\newcommand{\conv}{\operatorname{conv}}
\newcommand{\Span}{\operatorname{span}}
\newcommand{\argmin}{\operatorname{argmin}}
\newcommand{\ERMC}{{\hat f}_n^{\mathit{ERM}\mbox{{{\fontsize{7.6}{9.6}\selectfont{-}}}}C}}
\newcommand{\psinM}{\psi_n^{(C)}(M)}
\renewcommand{\hat}{\widehat}
\renewcommand{\epsilon}{\varepsilon}

\newcommand{\R}{\mathbb{R}}
\newcommand{\V}{\mathbb{V}}
\newcommand{\E}{\mathbb{E}}
\newcommand{\Pro}{\mathbb{P}}

\newcommand{\cA}{\mathcal{A}}
\newcommand{\cC}{\mathcal{C}}
\newcommand{\cD}{\mathcal{D}}
\newcommand{\cL}{\mathcal{L}}
\newcommand{\cX}{\mathcal{X}}
\newcommand{\cZ}{\mathcal{Z}}

\newtheorem{theorem}{Theorem}[section]
\newtheorem{proposition}[theorem]{Proposition}
\newremark{remark}[theorem]{Remark}

\makeatother

\begin{document}
\begin{frontmatter}

\title{Empirical risk minimization is optimal for the convex aggregation problem}
\runtitle{ERM in convex aggregation}

\begin{aug}
\author{\fnms{Guillaume} \snm{Lecu\'e}\corref{}\ead[label=e1]{guillaume.lecue@univ-mlv.fr}}
\runauthor{G. Lecu\'e} 
\address{CNRS, LAMA, Universit{\'e} Paris-Est Marne-la-vall{\'e}e,
77454 France.\\ \printead{e1}}
\end{aug}

\received{\smonth{10} \syear{2011}}
\revised{\smonth{1} \syear{2012}}

%
\begin{abstract}
Let $F$ be a finite model of cardinality $M$ and denote by $\conv(F)$
its convex hull.
The problem of convex aggregation is to construct a procedure having a
risk as close
as possible to the minimal risk over $\conv(F)$. Consider the bounded
regression model
with respect to the squared risk denoted by $R(\cdot)$. If $\ERMC$
denotes the empirical
risk minimization procedure over $\conv(F)$, then we prove that for
any $x>0$, with probability greater than $1-4\exp(-x)$,
\[
R(\ERMC)\leq\min_{f\in\conv(F)}R(f)+c_0\max \biggl(
\psi_n^{(C)}(M),\frac
{x}{n} \biggr),
\]
where $c_0>0$ is an absolute constant and $\psi_n^{(C)}(M)$ is the
optimal rate of convex aggregation
defined in (In \textit{Computational Learning Theory and Kernel Machines
(COLT-2003)} (2003) 303--313 Springer) by $\psinM=M/n$ when $M\leq
\sqrt {n}$ and $\psinM=\sqrt{\log (\mathrm{e}M/\sqrt{n} )/n}$ when
\mbox{$M>\sqrt{n}$}.
\end{abstract}

%

\end{frontmatter}

\section{Introduction and main results} \label{secintroduction}
Let $\cX$ be a probability space and let $(X,Y)$ and
$(X_1,Y_1),\ldots,(X_n,Y_n)$ be $n+1$ i.i.d. random variables with
values in
$\cX\times\R$. From the statistical point of view, the set
$\cD=\{(X_1,Y_1),\ldots,(X_n,Y_n)\}$ is the set of given data where the
$X_i$'s are usually considered as input data taking their values in
some space $\cX$ and the $Y_i$'s are some outputs or labels associated
with these inputs. We are interested in the prediction of $Y$
associated with a new observation~$X$. The data $\cD$ are thus used to
construct functions $f\dvtx \cX\rightarrow\R$ such that $f(X)$ provides a
good guess of $Y$. We measure the quality of this prediction by means
of the \textit{squared risk}
\[
R(f)=\E\bigl(Y-f(X)\bigr)^2,
\]
when $f$ is a real-valued function defined on $\cX$ and by
\[
R(\hat f)=\E \bigl[\bigl(Y-\hat f(X)\bigr)^2\vert\cD \bigr]
\]
when $\hat f$ is a function
constructed using the data $\cD$. For the sake of simplicity,
throughout this article, we restrict ourselves to functions
$f$ and random variables $(X,Y)$ for which
$|Y|\leq b$ and $|f(X)| \leq b$ almost surely, for some fixed $b\geq0$.
Note that $b$ does not have to be known from the statistician for the
construction of the procedures we are studying in this note.

Given a finite set $F$ of
real-valued measurable functions defined on $\cX$ (usually called a
\textit{dictionary}), there are three main types of aggregation problems:
\begin{enumerate}
\item \textit{Model selection aggregation}: construct a
procedure whose risk is as close as
possible to the risk of the best element in $F$ (cf. \cite
{Audibert1,MR2533466,MR2351101,MR2163920,DT07,MR1775638,MR1792783,MR2458184,MR2529440,TsyCOLT07,MR2051002,MR1790617,MR1762904}).
%
\item \textit{Convex aggregation}: construct a procedure whose risk is
as close as possible to the risk of the best function in the convex
hull of $F$ (cf. \cite
{MR2096215,MR2040405,MR2450298,MR2351101,MR1792783,TsyCOLT07,MR2044592,MR2356820}).
\item \textit{Linear aggregation}: construct a procedure whose risk is
as close as possible to the risk of the best function in the linear
span of $F$ (cf. \cite
{MR2351101,MR1775638,MR2329442,TsyCOLT07,AudCatRobust11}).
\end{enumerate}

In this note, we focus on the convex aggregation problem. We want to
construct a procedure $\tilde f$ for which, with high probability,
%
\begin{equation}
\label{eqOracle-ineq} R(\tilde f)\leq\min_{f\in\conv(F)}R(f)+
\psi_n(M),
\end{equation}
where $\psi_n(M)$ is called the residual term. The residual term is the
quantity that we want as small as possible. Results in expectation are
also of interest: construct a procedure $\tilde f$ such that $ \E
R(\tilde f)\leq\min_{f\in\conv(F)}R(f)+\psi_n(M)$.

In \cite{TsyCOLT07}, the author defined the \textit{optimal rates of
the convex aggregation}, by the smallest price in the minimax sense
that one has to pay to solve the convex aggregation problem. The
definition of \cite{TsyCOLT07} is given in expectation, as a function
of the cardinality $M$ of the dictionary $F$ and of the sample size
$n$. It has been proved in \cite{TsyCOLT07} (see also \cite{MR1792783}
and \cite{MR2044592}) that the optimal rate of convex aggregation is
\[
\psi_n^{(C)}(M)= \cases{ \displaystyle\frac{M}{n}, &\quad if $M\leq
\sqrt{n}$,
\cr
\displaystyle\sqrt{\frac{1}{n}\log \biggl(\frac{\mathrm{e}M}{\sqrt{n}} \biggr)}, &\quad if
$M>\sqrt{n}$. }
\]
This rate is defined up to some multiplying constant. Note that the
rate $ \psi_n^{(C)}(M)$ was achieved in \cite{TsyCOLT07} in expectation
for the Gaussian regression model with a known variance and a known
marginal distribution of the design. In \cite{MR2450298}, the authors
were able to remove these assumptions at a price of an extra $\log n$
factor for $1\leq M\leq\sqrt{n}$ (results are still in expectation).
Last year, there has been some striking results on different problems
of aggregation including the convex aggregation problem. To mention few
of them, we refer the reader to \cite{MR2816337,Rig11,WPGY12}. Finally, we also refer the reader to \cite
{MR2219712,MR2044592} for non-exact oracle inequalities (inequalities
like (\ref{eqOracle-ineq}) where $\min_{f\in\conv(F)}R(f)$ is
multiplied by a constant strictly larger than $1$) in the context of
convex aggregation.

A lower bound in deviation for the convex aggregation problem follows
from the arguments of \cite{TsyCOLT07}: there exist absolute positive
constants $c_0,c_1$ and $c_2$ such that for any sample cardinality
$n\geq1$, any cardinality of dictionary $M\geq1$ such that $\log
M\leq
c_0 n$, there exists a dictionary $F$ of size $M$ such that for any
aggregation procedure $\bar f_n$, there exists a random couple $(X,Y)$
such that $|Y|\leq b$ and $\max_{f\in F}|f(X)|\leq b$ a.s. and with
probability larger than $c_1$,
%
\begin{equation}
\label{eqlowerBound} R(\bar f_n)\geq\min_{f\in\conv(F) }R(f)+c_2
b^2\psi_n^{(C)}(M).
\end{equation}
This means that, from a minimax point of view, one cannot do better
than the rate $\psi_n^{(C)}(M)$ for the convex aggregation problem.
Therefore, any procedure achieving the rate $\psi_n^{(C)}(M)$ for any
dictionary $F$ and couple $(X,Y)$ such that $|Y|\leq b$ and $\max_{f\in
F}|f(X)|\leq b$ a.s. in an oracle inequality like (\ref{eqOracle-ineq})
is called an optimal procedure in deviation for the convex aggregation problem.

The procedure constructed in \cite{TsyCOLT07} achieves the rate $\psi_n^{(C)}(M)$ in
expectation (i.e., a~procedure satisfying (\ref
{eqOracle-ineq}) in expectation with the optimal residual term $\psi_n^{(C)}(M)$). An optimal procedure in deviation has been constructed
in Theorem 2.8.1 in \cite{HDR}. In both cases, the construction of
these optimal aggregation procedures require the aggregation of an
exponential number in $M$ of functions in $\conv(F)$ and thus cannot be
used in practice. On the other side, it would be much simpler and
natural to consider the classical procedure of\textit{ empirical risk
minimization} (cf. \cite{MR1641250}) over the convex hull of $F$ to
solve the convex aggregation problem:
%
\begin{equation}
\label{eqempirical-risk} \ERMC\in\argmin\limits_{f\in\conv(F)} R_n(f),
\qquad\mbox{where } R_n(f)=\frac
{1}{n}\sum
_{i=1}^n\bigl(Y_i-f(X_i)
\bigr)^2.
\end{equation}

In \cite{MR1792783,MR1775638,LM7}, the authors prove that, for every
$x>0$, with probability greater than $1-4\exp(-x)$
\[
R\bigl(\ERMC\bigr)\leq\min_{f\in\conv(F)}R(f)+ c_0\max \biggl(
\phi_n(M),\frac
{x}{n} \biggr), \quad \mbox{where }
\phi_n(M)=\min \biggl(\frac{M}{n},\sqrt {\frac
{\log M}{n}}
\biggr).
\]
The rate $\phi_n(M)$ behaves like the optimal rate $\psi_n^{(C)}(M)$
except for values of $M$ such that $n^{1/2} < M \leq c(\epsilon)
n^{1/2+\epsilon}$ for $\epsilon>0$ for which there is a logarithmic
gap. In this note, we were able to remove this logarithmic loss proving
that $\ERMC$ is indeed optimal for the convex aggregation problem.
Finally, note that in \cite{LM7}, the authors show that the rate $\psi_n^{(C)}(M)$ can be achieved by $\ERMC$ for any orthogonal dictionary
(i.e., such that $\forall f\neq g\in F, \E f(X) g(X)=0$). The
performance of ERM in the convex hull has been studied for an infinite
dictionary in \cite{MR2040405}. The resulting upper bounds, in the case
of a finite dictionary, is of the order of $M/n$ for every $n$ and~$M$.

Another motivation for this work comes from what is known about ERM in
the context of the three aggregation schemes mentioned above. It is
well known that ERM in $F$ is, in general, a suboptimal aggregation
procedure for the model selection aggregation problem (see \cite
{MR2458184,MR2451042} or \cite{LM2}). It is also known that ERM
in the linear span of $F$ is an optimal procedure for the linear
aggregation problem \cite{MR2329442} (cf. Theorem 13 and Example
1) or~\cite{AudCatRobust11}. Therefore, studying the performances of ERM in
the convex hull of $F$ in the context of convex aggregation can be seen
as an ``intermediate'' problem which remained open. In fact, a lot of
effort has been invested in finding any procedure that would be optimal
for the convex aggregation problem. For example, many boosting
algorithms (see \cite{MR1673273} or \cite{MR2420454} for recent results
on this topic) are based on finding the best convex combination in a
large dictionary (for instance, dictionaries consisting of ``decision
stumps''), while random forest algorithms can be seen as procedures
that try finding the best convex combination of decision trees. Thus,
finding an optimal procedure for the problem of convex aggregation for
a general dictionary is of high practical importance. In the following
result, we prove that empirical risk minimization is an optimal
procedure for the convex aggregation problem.

\renewcommand{\thetheorem}{A}
\begin{theorem}\label{theoA} There exists absolute
constants $c_0$ and $c_1$ such that the following holds. Let $F$ be a
finite dictionary of cardinality $M$ and $(X,Y)$ be a random couple of
$\cX\times\R$ such that $|Y|\leq b$ and $\max_{f \in F}|f(X)|\leq b$
a.s. for some $b>0$. Then, for any $x>0$, with probability greater than
$1-4\exp(- x)$
\[
R\bigl(\tilde{f}_n^{ERM-C}\bigr)\leq\min_{f\in{\conv}(F)}
R(f)+c_0b^2\max \biggl[\psi_n^{(C)}(M),
\frac{x}{n} \biggr].
\]
The optimality also holds in expectation:
\[
\E R\bigl(\tilde{f}_n^{ERM-C}\bigr)\leq\min_{f\in{\conv}(F)}
R(f)+c_1b^2\psi_n^{(C)}(M).
\]
\end{theorem}
%

\section{Preliminaries on isomorphic properties of functions
classes}\label{sectechnical-material}
We recall the machinery developed in \cite{MR2240689} to prove
isomorphic results between the empirical and actual structures of
functions classes.

Let $(\cZ,\sigma)$ be a measurable space, $Z,Z_1,\ldots,Z_n$ be $n+1$
i.i.d. random variables with values in $\cZ$ distributed according to
$P_Z$ and $G$ be a class of real-valued measurable functions defined on
$\cZ$. We consider the star shaped hull of $G$ in zero and its
localized set at some level $\lambda>0$:
\[
V(G)=\{\alpha g\dvt 0\leq\alpha\leq1, g\in G\} \quad\mbox{and}\quad V(G)_\lambda =
\bigl\{h\in V(G)\dvt P h\leq\lambda\bigr\}.
\]
For any functions class $H$ (in particular for $H$ being $G$, $V(G)$ or
$V(G)_\lambda$ for some $\lambda$), we denote $\llVert P-P_n\rrVert_H=\sup_{h\in H}|(P-P_n)h|$, where $Ph=\E h(Z)$ and $P_n h=n^{-1}\sum_{i=1}^nh(Z_i)$,
 $\sigma(H)=\sup_{h\in H}\sqrt{P h^2}$ and $\llVert
H\rrVert_\infty=\sup_{h\in H}\llVert h\rrVert_{L_\infty(P_Z)}$.
We also recall the
separability condition of \cite{MR2291502} (cf. Condition (M)) for
which Talagrand's concentration inequality holds:
\begin{enumerate}[(M)]
\item[(M)] There exists $G_0\subset G$ such that $G_0$ is countable and
for any $g\in G$, there exists a sequence $(g_k)_k$ in $G_0$ such that
for any $z\in\cZ$, $(g_k(z))_k$ tends to $g(z)$ when $k$ tends to infinity.
\end{enumerate}

\renewcommand{\thetheorem}{\arabic{section}.\arabic{theorem}}

\begin{theorem}[(\cite{MR2240689})]\label{theoshahar-peter}
There exists an absolute constant $c_0>0$ such that the following
holds. Let $G$ be a class of real-valued measurable functions defined
on $\cZ$ satisfying condition (\emph{M}) and such that $Pg^2\leq B Pg,\forall
g\in G$ for some constant $B>0$. Let $\lambda^*>0$ be such that
%
\begin{equation}
\label{eqfixed-point-general} \E\llVert P-P_n\rrVert_{V(G)_{\lambda^*}}
\leq(1/8)\lambda^*.
\end{equation}
For every $x>0$, with probability greater than $1-4\exp(-x)$, for every
$g\in G$,
\[
|Pg-P_ng|\leq(1/2)\max \bigl(Pg, \rho_n(x) \bigr),\qquad \mbox{where } \rho_n(x)=\max \biggl(\lambda^*,\frac{c_0(B +\llVert G\rrVert_\infty
)x}{n} \biggr).
\]
\end{theorem}
For the reader convenience, we recall the short proof of \cite{MR2240689}.

\begin{pf*}{Proof of Theorem \ref{theoshahar-peter}} Without loss of generality, we can assume that $G$ is
countable. From a limit argument, the result holds for classes of
functions satisfying condition (M).

Fix $\lambda>0$ and $x>0$, and note that by Talagrand's concentration
inequality (cf. \cite{MR1361756,MR1419006,MR2319879,MR2135312,MR1890640}), with probability
larger than
$1-4\exp(-x)$,
%
\begin{equation}
\label{eqAdamcjak} \llVert P - P_n\rrVert_{V(G)_\lambda} \leq2 \E
\llVert P-P_n\rrVert_{V(G)_{\lambda}} + K \sigma\bigl(V(G)_{\lambda}
\bigr) \sqrt{\frac x n} + K \bigl\llVert V(G)_{\lambda}\bigr
\rrVert_\infty\frac x n,
\end{equation}
where $K$ is an absolute constant.
Clearly, we have $\llVert V(G)_\lambda\rrVert_\infty\leq\llVert
G\rrVert_\infty$ and
\[
\sigma^2\bigl(V(G)_{\lambda}\bigr) = \sup \bigl(P(\alpha
g)^2 \dvt  0 \leq\alpha \leq1, g \in G, P( \alpha g) \leq\lambda \bigr)
\leq B\lambda.
\]

Moreover, since $V(G)$ is star-shaped, $\lambda>0\rightarrow\phi
(\lambda)= \E\llVert P-P_n\rrVert_{V(G)_\lambda}/\lambda$ is
non-increasing, and since $\phi(\lambda^*)\leq1/8$ and
$\rho_n(x)\geq\lambda^*$ then
\[
\E\llVert P-P_n\rrVert_{V(G)_{\rho_n(x)}}\leq(1/8)\rho_n(x).
\]
Combined with \eqref{eqAdamcjak}, there
exists an event $\Omega_0(x)$ of probability greater than
$1-4\exp(-x)$, and on $\Omega_0(x)$,
\[
\llVert P-P_n\rrVert_{V(G)_{\rho_n(x)}} \leq (1/4)\rho_n(x)
+ K \sqrt{\frac{B\rho_n(x) x}{n}} + K\frac{\llVert G\rrVert_\infty
x}{n} \leq (1/2)
\rho_n(x)
\]
as long as $c_0\geq64(K^2+K)$.
Hence, on $\Omega_0(x)$, if $g\in V(G)$ satisfies that $P g\leq\rho_n(x)$, then $|Pg - P_ng|\leq(1/2)\rho_n(x)$. Moreover, if $g\in V(G)$
is such that
$Pg>\rho_n(x)$, then
$h=\rho_n(x)g/P g \in V(G)_{\rho_n(x)}$; hence
$|Ph-P_nh|\leq
(1/2)\rho_n(x)$, and so in both cases $|Pg-P_ng|\leq(1/2)\max
(Pg,\rho_n(x) )$.
\end{pf*}

Therefore, if one applies Theorem \ref{theoshahar-peter} to obtain
isomorphic properties between the empirical and actual structures, one
has to check the condition $Pg^2\leq B Pg,\forall g\in G$, called the
Bernstein condition in \cite{MR2240689}, and to find a point $\lambda^*$ satisfying (\ref{eqfixed-point-general}).

A point $\lambda^*$ such that (\ref{eqfixed-point-general}) holds can
be found thanks to the peeling argument of \cite{BMN}: for any
$\lambda>0$,
\[
V(G)_\lambda\subset\bigcup_{i=0}^\infty
\bigl\{\alpha g\dvt 0\leq\alpha \leq2^{-i}, g\in G, Pg\leq2^{i+1}
\lambda \bigr\}
\]
which implies
%
\begin{equation}
\label{eqpeeling-shahar} \E\llVert P-P_n\rrVert_{V(G)_\lambda}\leq
\sum_{i=0}^\infty 2^{-i}\E\llVert
P-P_n\rrVert_{G_{2^{i+1}\lambda}},
\end{equation}
where, for any $\mu>0$, $G_\mu=\{g\in G\dvt Pg\leq\mu\}$. Then if
$\lambda^*>0$ is such that $\lambda^*/8$ upper bounds the RHS in (\ref
{eqpeeling-shahar}) this point also satisfies (\ref{eqfixed-point-general}).

The Bernstein condition usually follows from some convexity argument.
For instance, it is now standard to check the Bernstein condition for
the excess loss functions class $\cL_F=\{\cL_f\dvt f\in F\}$ associated
with a convex model $F$ with respect to the squared loss function $\ell_f(x,y)=(y-f(x))^2,\forall x\in\cX,y\in\R$,
 where $f^*_F\in\argmin_{f\in F}\E(Y-f(X))^2$ and $\cL_f=\ell_f-\ell_{f^*_F}$. Indeed, if $F$
is a convex set of functions and $(X,Y)$ is a random couple on $\cX
\times\R$ such that $|Y|\leq b$ and $\sup_{f\in F}|f(X)|\leq b$ a.s.
then it follows from convexity and definition of $f^*_F$ that for any
$f\in F$, $\E [(f^*_F(X)-Y)(f^*_F(X)-f(X)) ]\leq0$ and so
%
\begin{equation}
\label{eqone-side-Bernsetin-1} \E\cL_f=2\E\bigl(f^*_F(X)-f(X)
\bigr) \bigl(Y-f^*_{F}(X)\bigr)+\E\bigl(f^*_{F}(X)-f(X)
\bigr)^2\geq \E \bigl(f^*_F(X)-f(X)\bigr)^2.
\quad\end{equation}
Moreover, since $|Y|\leq b$ and $\sup_{f\in F}|f(X)|\leq b$ a.s. then
%
\begin{equation}
\label{one-side-Bernstein-2} \E\cL_f^2=\E
\bigl(2Y-f^*_{F}(X)-f(X) \bigr)^2 \bigl(f(X)-f^*_{F}(X)
\bigr)^2\leq(4b)^2\E \bigl(f(X)-f^*_{F}(X)
\bigr)^2.
\end{equation}
Therefore, any $f$ in $F$ is such that $\E\cL_f^2\leq(4b)^2\E\cL_f$.

\section{\texorpdfstring{Proof of Theorem \protect\ref{theoA}}{Proof of Theorem A}}\label{secUpper-bound-MS-aggregation}

The proof of Theorem \ref{theoA} for the case $M\leq\sqrt{n}$ is now very
classical and can be found in \cite{MR2329442} (cf. Theorem 13 and
Example 1). Nevertheless, we reproduce here this short proof in order
to provide a self-contained note. The proof for the case $M>\sqrt{n}$
is more tricky and relies on isomorphic properties of an exponential
number of segments in $\conv(F)$ together with Maurey's empirical
method (cf. \cite{MR659306,MR810669}) which was first used in the
context of convex aggregation in \cite{MR1775638} and \cite{TsyCOLT07}.
Note that segments are models of particular interest in Learning theory
because they are convex models (in particular, they satisfy the
Bernstein condition) and they are of small complexity (essentially the
same complexity as a model of cardinality two). On the contrary to the
classical entropy based approach which essentially consists in
approximating a set by finite sets, approaching models by union of
segments may be of particular interest in Learning theory beyond the
convex aggregation problem. Note that finite models have no particular
geometrical structure and therefore are somehow ``bad models'' as far
as ERM procedures are concerned.

Proofs are given for the deviation result of Theorem \ref{theoA}. The result in
expectation of Theorem~\ref{theoA} follows from a direct integration argument.

\subsection{\texorpdfstring{The case $M>\sqrt{n}$}{The case M > square root of n}}\label{secMgeqSquaren}

We apply Theorem \ref{theoshahar-peter} to excess loss functions
classes indexed by segments. First, note that segments of bounded
functions are functions classes satisfying condition (M). We consider a
set $\cC^\prime=\{g_1,\ldots,g_N\}$ of real-valued measurable functions
defined on $\cX$ such that $\max_{g\in\cC^\prime}|g(X)|\leq b$ a.s.
For every $i,j\in\{1,\ldots,N\}$, we consider the segment
$[g_i,g_j]=\{
\theta g_i+(1-\theta)g_j\dvt 0\leq\theta\leq1\}$ and take $g^*_{ij}\in
\argmin_{g\in[g_i,g_j]} R(g)$ where $R(\cdot)$ is the squared risk. We
consider the excess loss functions class
\[
\cL^{ij}=\bigl\{\cL_g^{ij}\dvt g
\in[g_i,g_j]\bigr\},\qquad \mbox{where } \cL_g^{ij}=
\ell_g-\ell_{g^*_{ij}}
\]
for $\ell_g(x,y)=(y-g(x))^2,\forall x\in\cX,y\in\R$.

As a consequence of convexity of segments, we have for any $g\in
[g_i,g_j]$, $\E (\cL^{ij}_g )^2\leq(4b)^2\E\cL^{ij}_g$ (cf.
(\ref{eqone-side-Bernsetin-1}) and (\ref{one-side-Bernstein-2}) in
Section \ref{sectechnical-material}). This implies that the functions
class $\cL^{ij}$ satisfies the Bernstein condition of Theorem \ref
{theoshahar-peter}. Now, it remains to find $\lambda^*>0$ such that
$\E
\llVert P-P_n\rrVert_{V(\cL^{ij})_{\lambda^*}}\leq(1/8)\lambda^*$. Let
$\mu>0$
and $\epsilon_1,\ldots,\epsilon_n$ be $n$ i.i.d. Rademacher variables.
Note that for any $g\in[g_i,g_j]$, $P\cL^{ij}_g\geq
P(g-g_{ij}^*)^2=\E
(g(X)-g_{ij}^*(X) )^2$ (cf. (\ref{eqone-side-Bernsetin-1})). It
follows from the symmetrization argument and the contraction principle
(cf. \cite{LT91}, page 95) that if $g_i\neq g_j$ then
\begin{eqnarray*}
\E\llVert P-P_n\rrVert_{(\cL^{ij})_\mu}&\leq&2\E\sup_{g\in
[g_i,g_j]\dvt P\cL
^{ij}_g\leq\mu} \Biggl|
\frac{1}{n}\sum_{k=1}^n
\epsilon_k\cL^{ij}_g(X_k,Y_k)
\Biggr|
\\
&\leq& 8b\E\sup_{g\in[g_i,g_j]\dvt P\cL^{ij}_g\leq\mu} \Biggl|\frac
{1}{n}\sum
_{k=1}^n\epsilon_k\bigl(g(X_k)-g^*_{ij}(X_k)
\bigr) \Biggr|
\\
&\leq& 8b\E\sup_{g\in[g_i,g_j]\dvt P(g-g^*_{ij})^2\leq\mu} \Biggl|\frac
{1}{n}\sum
_{k=1}^n\epsilon_k\bigl(g(X_k)-g^*_{ij}(X_k)
\bigr) \Biggr|
\\
&=& 8b\E\sup_{g\in[g_i,g_j]-g^*_{ij}\dvt Pg^2\leq\mu} \Biggl|\frac
{1}{n}\sum_{k=1}^n
\epsilon_kg(X_k) \Biggr|\\
&\leq& 8b\E\sup_{g\in\Span
(g_i-g_j)\dvt Pg^2\leq\mu} \Biggl|
\frac{1}{n}\sum_{k=1}^n
\epsilon_kg(X_k)\Biggr |
\\
&=& \frac{8b\sqrt{\mu}}{P(g_i-g_j)^2}\E \Biggl|\frac{1}{n}\sum_{k=1}^n
\epsilon_k (g_i-g_j) (X_k) \Biggr|\\
&\leq&\frac{8b\sqrt{\mu
}}{P(g_i-g_j)^2} \Biggl(\E \Biggl(\frac{1}{n}\sum
_{k=1}^n\epsilon_k (g_i-g_j)
(X_k) \Biggr)^2 \Biggr)^{1/2}
\\
&=& 8b\sqrt{\frac{\mu}{n}}.
\end{eqnarray*}
Note that when $g_i=g_j$ the result is also true. Now, we use the
peeling argument of (\ref{eqpeeling-shahar}) to obtain
\[
\E\llVert P-P_n\rrVert_{V(\cL^{ij})_\lambda}\leq\sum
_{k=0}^\infty 2^{-k}\E \llVert
P-P_n\rrVert_{(\cL^{ij})_{2^{k+1}\lambda}}\leq\sum_{k=0}^\infty
2^{-k}8b \sqrt{\frac{2^{k+1}\lambda}{n}}\leq c_0 b \sqrt{
\lambda/n}.
\]
Therefore, for $\lambda^*=(8c_0b)^2/n$, we have $ \E\llVert
P-P_n\rrVert_{V(\cL
^{ij})_{\lambda^*}}\leq(1/8)\lambda^*$.

Now, we can apply Theorem \ref{theoshahar-peter} to the family of
excess loss functions classes $(\cL^{ij})_{1\leq i,j\leq N}$ together
with a union bound to obtain the following result.

\begin{proposition}\label{propUniform-Isomoprhism-segment}There exists
an absolute constant $c_0>0$ such that the following holds. Let $\cC^\prime=\{g_1,\ldots,g_N\}$ be a set of measurable real-valued
functions defined on $\cX$. Let $(X,Y)$ be a random couple with values
in $\cX\times\R$ such that $|Y|\leq b$ and $\max_{g\in\cC^\prime
}|g(X)|\leq b$ a.s. For any $x>0$, with probability greater than
$1-4\exp(-x)$, for any $i,j\in\{1,\ldots,N\}$ and any $g\in[g_i,g_j]$,
\[
\bigl|P\cL_g^{ij}-P_n\cL^{ij}_g
\bigr|\leq(1/2)\max \bigl(P\cL_g^{ij},\gamma(x) \bigr),\qquad \mbox{where } \gamma(x)=\frac{c_0
b^2(x+2\log N)}{n}.
\]
\end{proposition}

Now, we want to apply the isomorphic result of Proposition \ref
{propUniform-Isomoprhism-segment} to a wisely chosen subset $\cC^\prime
$ of $\cC=\conv(F)$. For that, we consider the integer
\[
m= \biggl\lceil\sqrt{\frac{n}{\log (\mathrm{e}M/\sqrt{n} )}} \biggr\rceil
\]
and the set $\cC^\prime$ is defined by
\[
\cC^\prime= \Biggl\{\frac{1}{m}\sum_{i=1}^mh_i\dvt h_1,
\ldots,h_m\in F \Biggr\}.
\]
The set $\cC^\prime$ is an approximating set of the convex hull
$\conv
(F)$. We will, for instance, use the following approximation property:
%
\begin{equation}
\label{eqMaurey-1} \min_{f\in\cC^\prime}R(f)\leq\min_{f\in\cC}R(f)+
\frac{4b^2}{m}.
\end{equation}
Indeed, to obtain such a result, we use Maurey's empirical method. Let
$f^*_\cC\in\argmin_{f\in\cC}R(f)$ and denote $f^*_\cC=\sum_{j=1}^M\lambda_j f_j$ where $\lambda_j\ge0,\forall j=1,\ldots,M$ and
$\sum_{j=1}^M\lambda_j=1$. Consider a random variable $\Theta\dvtx \Omega
\rightarrow F$ such that $\Pro[\Theta=f_j]=\lambda_j,\forall
j=1,\ldots
,M$ and let $\Theta_1,\ldots,\Theta_m$ be $m$ i.i.d. random variables
distributed according to $\Theta$ and independent of $(X,Y)$. Denote by
$\E_\Theta$ the expectation with respect to $\Theta_1,\ldots,\Theta_m$.
Since $\E_\Theta\Theta_j=f^*_\cC$ for any $j=1,\ldots,m$, we have
\begin{eqnarray*}
\min_{f\in\cC^\prime}R(f)&\leq&\E_\Theta R \Biggl(\frac{1}{m}\sum
_{j=1}^m\Theta_j \Biggr)=
\E_\Theta\E \Biggl(\frac{1}{m}\sum_{j=1}^m
\Theta_j(X)-Y \Biggr)^2
\\
&=&\E \Biggl(\frac{1}{m^2}\sum_{j,k=1}^m
\E_\Theta\bigl(Y-\Theta_j(X)\bigr) \bigl(Y-
\Theta_k(X)\bigr) \Biggr)=R\bigl(f^*_\cC\bigr)+
\frac{\E\V_\Theta(Y-\Theta(X))}{m},
\end{eqnarray*}
where $\V_\Theta$ stands for the variance symbol with respect to
$\Theta
$. Equation (\ref{eqMaurey-1}) follows since $|Y|\leq b$ and $\max_{f\in F}|f(X)|\leq b$ a.s.

Denote by $N=|\cC^\prime|$ the cardinality of $\cC^\prime$ and by
$g_1,\ldots,g_N$ the functions in $\cC^\prime$. For simplicity, assume
that $R(g_1)=\min_{g\in\cC^\prime}R(g)$. Thanks to \cite{MR810669} for
the first inequality and \cite{DGL96}, page 218, or \cite{MR2319879},
Proposition 2, for the second inequality, we know that
%
\begin{equation}
\label{eqCarl-denombrement} \bigl|\cC^\prime\bigr|=N\leq %
\pmatrix{
M+m-1
\cr
m }\leq \biggl(\frac{2\mathrm{e}M}{m} \biggr)^m.
\end{equation}

Let $x>0$. Consider the event $\Omega(x)\subset\Omega$ such that the
following isomorphic property holds for all the segments $[g_1,g_j],
j=1,\ldots,N$:
%
\begin{equation}
\label{eqiso-property-all-segments} \bigl|P_n\cL^{1j}_g-P
\cL^{1j}_g \bigr|\leq(1/2)\max \bigl(P\cL^{1j}_g,
\gamma(x) \bigr)\qquad \forall g\in[g_1,g_j],
\end{equation}
where we recall that $\cL^{1j}_g=\ell_g-\ell_{g^*_{1j}}$ is the excess
loss function of $g\in[g_1,g_j]$ for the model $[g_1,g_j]$ and
\[
\gamma(x)=\frac{c_0b^2(x+2\log N)}{n}.
\]
Thanks to Proposition \ref{propUniform-Isomoprhism-segment}, we know
that $\Pro[\Omega(x)]\geq1-4\exp(-x)$.

We are going to work on the event $\Omega(x)$ but for the moment, we
use a second time Maurey's empirical method. Fix $X_1,\ldots,X_n$ and
write $\ERMC=\sum_{j=1}^M\beta_j f_j$. Consider a random variable
$\Theta\dvtx  \Omega^\prime\rightarrow F$ defined on an other probability
space $(\Omega^\prime,\cA^\prime,\Pro^\prime)$ such that $\Pro^\prime
[\Theta=f_j]=\beta_j,\forall j=1,\ldots,M$ and let $\Theta_1,\ldots
,\Theta_m$ be $m$ i.i.d. random variables having the same probability
distribution as $\Theta$. Once again, denote by $\E_\Theta^\prime$ the
expectation with respect to $\Theta_1,\ldots,\Theta_m$ and by $\V_\Theta
$ the variance with respect to $\Theta$. Since $\E_\Theta^\prime
\Theta_j=\ERMC$ for any $j=1,\ldots,m$, it follows from the same method used
to obtain (\ref{eqMaurey-1}) that
%
\begin{equation}
\label{eqMaurey-2} \E_\Theta^\prime R \Biggl(\frac{1}{m}
\sum_{j=1}^m\Theta_j
\Biggr)=R\bigl(\ERMC \bigr)+\frac{\E\V_\Theta^\prime(Y-\Theta(X))}{m}
\end{equation}
and the same holds for the empirical risk:
%
\begin{equation}
\label{eqMaurey-3} \E_\Theta^\prime R_n \Biggl(
\frac{1}{m}\sum_{j=1}^m
\Theta_j \Biggr)=R_n\bigl(\ERMC \bigr)+\frac{1}{m} \Biggl(
\frac{1}{n}\sum_{i=1}^n
\V_\Theta^\prime \bigl(Y_i-\Theta (X_i)
\bigr) \Biggr).
\end{equation}
Consider the following notation:
\[
g_\Theta=\frac{1}{m}\sum_{j=1}^m
\Theta_j \quad\mbox{and}\quad i_\Theta\in \{ 1,\ldots,N\} \mbox{ such that } g_{i_\Theta}=g_\Theta.
\]
Note that $g_\Theta$ is a random point in $\cC^\prime$ (as a measurable
function from $\Omega^\prime$ to $\cC^\prime$) and that, on the event
$\Omega(x)$, the following isomorphic property on the segment
$[g_1,g_\Theta]$ holds:
%
\begin{equation}
\label{eqiso-property-segments-Theta} \bigl|P_n\cL^{1i_\Theta}_g-P
\cL^{1i_\Theta}_g \bigr|\leq(1/2)\max \bigl(P\cL^{1i_\Theta}_g,
\gamma(x) \bigr)\qquad \forall g\in[g_1,g_{i_\Theta}].
\end{equation}

First note that for every $\Theta_1,\ldots,\Theta_m$, we have
%
\begin{equation}
\label{eqequa-1} R\bigl(\ERMC\bigr)=R\bigl(g^*_{1i_\Theta}\bigr)+R(g_\Theta)-R
\bigl(g^*_{1i_\Theta}\bigr)+R\bigl(\ERMC \bigr)-R(g_\Theta).
\end{equation}
By definition of $g^*_{1i_\Theta}\in\argmin_{g\in[g_{i_\Theta},g_1]}
R(g)$, we have $R(g^*_{1i_\Theta})\leq R(g_1)=\min_{g\in\cC^\prime
}R(g)$ and according to (\ref{eqMaurey-1}), we have $\min_{f\in\cC
^\prime}R(f)\leq\min_{f\in\cC}R(f)+(4b^2)/m$. Therefore, it follows
from (\ref{eqequa-1}) that
%
\begin{equation}
\label{eqequa-2} R\bigl(\ERMC\bigr)\leq\min_{f\in\cC}R(f)+\frac{4b^2}{m}+P
\cL^{1 i_\Theta
}_{g_\Theta}+R\bigl(\ERMC\bigr)-R(g_\Theta).
\end{equation}
On the event $\Omega(x)$, we use (\ref{eqiso-property-segments-Theta})
to obtain for every $\Theta_1,\ldots,\Theta_m$
\[
R\bigl(\ERMC\bigr)\leq\min_{f\in\cC}R(f)+\frac{4b^2}{m}+2P_n
\cL^{1
i_\Theta
}_{g_\Theta}+\gamma(x)+R\bigl(\ERMC\bigr)-R(g_\Theta).
\]
Moreover, by definition of $\ERMC$, we have
\[
P_n\cL^{1 i_\Theta}_{g_\Theta}=R_n(g_\Theta)-R_n
\bigl(g_{1i_\Theta
}^*\bigr)\leq R_n(g_\Theta)-R_n\bigl(
\ERMC\bigr).
\]
Therefore, on the event $\Omega(x)$, we have for every $\Theta_1,\ldots
,\Theta_m$
\begin{eqnarray*}
R\bigl(\ERMC\bigr)&\leq&\min_{f\in\cC}R(f)+\frac{4b^2}{m}+\gamma(x)
\\
&&{}+2 \bigl(R_n(g_\Theta)-R_n\bigl(\ERMC\bigr) \bigr)+R\bigl(
\ERMC\bigr)-R(g_\Theta).
\end{eqnarray*}
In particular, one can take the expectation with respect to $\Theta_1,\ldots,\Theta_m$ (defined on $\Omega^\prime$) in the last
inequality. We have on $\Omega(x)$,
\begin{eqnarray*}
R\bigl(\ERMC\bigr)&\leq&\min_{f\in\cC}R(f)+\frac{4b^2}{m}+\gamma(x)
\\
&&{}+2\E_\Theta^\prime \bigl(R_n(g_\Theta)-R_n\bigl(
\ERMC\bigr) \bigr)+\E_\Theta^\prime \bigl(R\bigl(\ERMC\bigr)-R(g_\Theta)
\bigr).
\end{eqnarray*}
Thanks to (\ref{eqMaurey-2}), we have $\E_\Theta^\prime
(R(\ERMC
)-R(g_\Theta) )\leq0$ and it follows from (\ref{eqMaurey-3})
that $\E_\Theta^\prime (R_n(g_\Theta)-R_n(\ERMC) )\leq(2b)^2/m$.
Therefore, on the event $\Omega(x)$, we have
\[
R\bigl(\ERMC\bigr)\leq\min_{f\in\cC}R(f)+\frac{8b^2}{m}+\gamma(x)\leq
\min_{f\in
\cC}R(f)+c_1 b^2\max \biggl(\psinM,
\frac{x}{n} \biggr),
\]
where the last inequality follows from (\ref{eqCarl-denombrement}) and
the definition of $m$.

\subsection{\texorpdfstring{The case $M\leq\sqrt{n}$}
{The case M <= square root of n$}}\label{secMleqSquaren}\enlargethispage{3pt}

We use the strategy developed in \cite{MR2240689} together with the one
of \cite{MR2329442} (cf. Example 1) to prove Theorem \ref{theoA} in the case
$M\leq\sqrt{n}$. Define $\cC=\conv(F)$ and $\cL_{\cC}=\{\cL_f\dvt f\in\cC\}
$ the excess loss class associated with $\cC$ where $\cL_f=\ell_f-\ell_{f^*_{\cC}},\forall f\in\cC$ and $f^*_\cC\in\argmin_{f\in\cC}R(f)$.

Let $x>0$. Assume that we can find some $\rho_n(x)>0$ such that with
probability greater than $1-4\exp(-x)$, for any $f\in\cC$,
%
\begin{equation}
\label{eqiso-convex-hull} |P_n\cL_f-P\cL_f
|\leq(1/2)\max \bigl(P\cL_f,\rho_n(x) \bigr).
\end{equation}
Then, the ERM over $\conv(F)$ would satisfy with probability greater
than $1-4\exp(-x)$,
\[
R\bigl(\ERMC\bigr)-\min_{f\in\conv(F)}R(f)=P\cL_{\ERMC}\leq2P_n
\cL_{\ERMC
}+\rho_n(x)\leq\rho_n(x).
\]
This means that if we can prove some isomorphic properties between the
empirical and the actual structures of the functions class $\cL_\cC$
like in (\ref{eqiso-convex-hull}), then we can derive oracle
inequalities for $\ERMC$. This is the strategy used in \cite{MR2240689}
that we follow here.

According to Theorem \ref{theoshahar-peter} in Section \ref
{sectechnical-material}, a function $\rho_n(x)$ satisfying (\ref
{eqiso-convex-hull}) can be constructed if we prove that $\cL_\cC$
satisfies some Bernstein condition and if we find some fixed point
$\lambda^*>0$ such that $\E\llVert P-P_n\rrVert_{V(\cL_\cC
)_{\lambda
^*}}\leq
(1/8)\lambda^*$. The Bernstein condition follows from the convexity of
$\conv(F)$ and the strategy used in Section \ref
{sectechnical-material}: for any $f\in\cC, P\cL_f^2\leq(4b)^2 P\cL_f$.

We use the peeling argument of Section \ref{sectechnical-material}
together with the following observations due to \cite{MR2329442} (cf.
Example 1) to find a fixed point $\lambda^*$. Let $S$ be the linear
subspace of $L^2(P_X)$ spanned by the dictionary $F$ and consider an
orthonormal basis $(e_1,\ldots,e_{M^\prime})$ of $S$ in $L^2(P_X)$
(where $M^\prime={\dim}(S)\leq M$). For any $\mu>0$, it follows from
the symmetrization argument and the contraction principle (cf. Chapter
4 in \cite{LT91}) that
\begin{eqnarray*}
\E\llVert P-P_n\rrVert_{(\cL_\cC)_\mu}&\leq&8b\E\sup_{f\in
S\dvt Pf^2\leq
\mu} \Biggl|
\frac{1}{n}\sum_{i=1}^n
\epsilon_i f(X_i)\Biggr |
\\[-2pt]
&\leq&8b\E\sup_{\beta\in\R^{M^\prime}\dvt \|\beta\|_{2} \leq\sqrt
{\mu}} \Biggl|\frac{1}{n}\sum
_{i=1}^n\epsilon_i \Biggl(\sum
_{j=1}^{M^\prime
}\beta_je_j(X_i)
\Biggr) \Biggr|
\\[-2pt]
&\leq&8b\sqrt{\mu} \E \Biggl(\sum_{j=1}^{M^\prime}
\Biggl(\frac
{1}{n}\sum_{i=1}^n
\epsilon_ie_j(X_i) \Biggr)^2
\Biggr)^{1/2} \leq8b\sqrt{\frac{M^\prime\mu}{n}}.
\end{eqnarray*}
We use the peeling argument of (\ref{eqpeeling-shahar}) to prove that
for $\lambda^*=c_0b^2M/n$ and $c_0$ an absolute constant large enough,
we have indeed $\E\llVert P-P_n\rrVert_{V(\cL_\cC)_{\lambda
^*}}\leq
(1/8)\lambda^*$.

Now, it follows from Theorem \ref{theoshahar-peter} that for any $x>0$,
with probability greater than $1-4\exp(-x)$,
\[
R\bigl(\ERMC\bigr)\leq\min_{f\in\cC}R(f)+c_1b^2\max \biggl(
\frac{M}{n},\frac
{x}{n} \biggr).
\]
This concludes the proof for the case $M\leq\sqrt{n}$.\eject

\begin{remark} We did not use the condition $M\leq\sqrt{n}$ in the
last proof. In fact, the result holds in the following more general
framework. Let $\Lambda$ be any closed convex subset of $\R^M$ and for
any dictionary $F=\{f_1,\ldots,f_M\}$ denote by $\Lambda(F)$ the set of
all functions $\sum_{j=1}^M\lambda_j f_j$ when $(\lambda_1,\ldots
,\lambda_M)^\top\in\Lambda$. Let $(X,Y)$ be a random couple with values
in $\cX\times\R$ such that $|Y|\leq b$ and $\max_{f\in F}|f(X)|\leq b$
a.s. Consider the ERM procedure
\[
\hat f_n\in\argmin\limits_{f\in\Lambda(F)}R_n(f).
\]
Then, it follows from Theorem \ref{theoshahar-peter} and the argument
used previously in this section that for any $x>0$, with probability
greater than $1-4\exp(-x)$,
\[
R(\hat f_n)\leq\min_{f\in\Lambda(F)}R(f)+c_1b^2
\max \biggl(\frac
{M}{n},\frac{x}{n} \biggr).
\]
The same result can be found in \cite{AudCatRobust11} under very weak
moment assumptions.
\end{remark}

\section*{Acknowledgements}
We would like to thank Alexandre Tsybakov for helping us for the
presentation of this result.
Supported by French Agence Nationale de la Recherche ANR Grant
\textsc{``Prognostic''} ANR-09-JCJC-0101-01.



\printhistory

\end{document}